\documentclass[a4, 10pt]{amsart}
\usepackage{amssymb}
\usepackage{amstext}
\usepackage{amsmath}
\usepackage{amscd}
\usepackage{latexsym}
\usepackage{amsfonts}

\theoremstyle{plain}
\newtheorem{thm}{Theorem}[section]
\newtheorem*{thm*}{Theorem}
\newtheorem*{cor*}{Corollary}

\newtheorem{prop}[thm]{Proposition}
\newtheorem{lem}[thm]{Lemma}

\newtheorem{cor}[thm]{Corollary}
\newtheorem{claim}{Claim}
\newtheorem*{claim*}{Claim}

\theoremstyle{definition}
\newtheorem{defn}[thm]{Definition}
\newtheorem{ex}[thm]{Example}
\newtheorem{rem}[thm]{Remark}
\newtheorem{conj}[thm]{Conjecture}

\theoremstyle{remark}
\newtheorem*{pf}{{\sl Proof}}

\newtheorem*{cpf}{{\sl Proof of Claim}}

\numberwithin{equation}{thm}
\def\Hom{\mathrm{Hom}}

\def\Ext{\mathrm{Ext}}

\def\Tor{\mathrm{Tor}}

\def\mod{\mathrm{mod}}

\def\m{\mathfrak m}
\def\n{\mathfrak n}

\def\p{\mathfrak p}
\def\q{\mathfrak q}

\def\P{\mathfrak P}

\def\N{\Bbb N}

\def\C{\Bbb C}

\def\Supp{\mathrm{Supp}}
\def\Ann{\mathrm{Ann}}
\def\Ass{\mathrm{Ass}}

\def\Min{\mathrm{Min}}
\def\pd{\mathrm{pd}}

\def\grade{\mathrm{grade}}

\def\Spec{\mathrm{Spec}}
\def\Sing{\mathrm{Sing}}

\def\X{{\mathcal X}}

\def\G{{\mathcal G}}
\def\M{\mathfrak M}

\def\ind{\mathrm{ind}}
\def\CM{\mathrm{CM}}
\def\NF{\mathrm{NF}}

\begin{document}

\title[Indecomposable totally reflexive modules]{An uncountably infinite number of indecomposable totally reflexive modules}
\author{Ryo Takahashi}
\address{Department of Mathematics, School of Science and Technology, Meiji University, Kawasaki 214-8571, Japan}
\email{takahasi@math.meiji.ac.jp}
\keywords{Cohen-Macaulay ring, Countable Cohen-Macaulay type, Totally reflexive, Semidualizing}
\subjclass[2000]{13C14, 16G60}
\begin{abstract}
A few years ago, Huneke and Leuschke proved a theorem which solved a conjecture of Schreyer.
It asserts that an excellent Cohen-Macaulay local ring of countable Cohen-Macaulay type which is complete or has uncountable residue field has at most a one-dimensional singular locus.
In this paper, it is verified that the assumption of the excellent property can be removed, and the theorem is considered over an arbitrary local ring.
The main purpose of this paper is to prove that the existence of a certain prime ideal and a certain totally reflexive module implies the existence of an uncountably infinite number of isomorphism classes of indecomposable totally reflexive modules.
\end{abstract}
\maketitle
\section{Introduction}

Throughout the present paper, we assume that all rings are commutative and noetherian, and that all modules are finitely generated.

A Cohen-Macaulay local ring is said to be of finite (resp. countable) Cohen-Macaulay type if there exist only finitely (resp. countably) many isomorphism classes of indecomposable maximal Cohen-Macaulay modules.
The property of finite Cohen-Macaulay type has been deeply studied for almost thirty years.
The following theorem is a well-known result concerning this property; it was proved by Auslander \cite{Auslander2} in the complete case, Leuschke and Wiegand \cite{LW} in the excellent case, and Huneke and Leuschke \cite{HL} in the general case:

\begin{thm}[Auslander-Huneke-Leuschke-Wiegand]\label{}
Let $R$ be a Cohen-Macaulay local ring of finite Cohen-Macaulay type.
Then $R$ has at most an isolated singularity.
\end{thm}

It is obvious that a local ring $R$ has an isolated singularity if and only if $\dim R/\p =0$ for any $\p\in\Sing\,R$, where $\Sing\,R$ denotes the singular locus of $R$.
Schreyer \cite{Schreyer} conjectured that the singular locus of a Cohen-Macaulay local ring of countable Cohen-Macaulay type consists of prime ideals of lower dimension:

\begin{conj}[Schreyer]\label{}
Let $R$ be an analytic Cohen-Macaulay local $\C$-algebra of countable Cohen-Macaulay type.
Then $R$ has at most a one-dimensional singular locus, namely, $\dim R/\p\leq 1$ for any $\p\in\Sing\,R$.
\end{conj}

This conjecture was recently proved by Huneke and Leuschke \cite{HL2}; they actually proved a stronger statement than Schreyer's conjecture: 

\begin{thm}[Huneke-Leuschke]\label{hl}
Let $(R, \m , k)$ be an excellent Cohen-Macaulay local ring of countable Cohen-Macaulay type.
Suppose either that $R$ is complete or that $k$ is uncountable.
Then $\dim R/\p\leq 1$ for any $\p\in\Sing\,R$.
\end{thm}

Actually, in the above theorem, the assumption that $R$ is excellent is not necessary; we will show that in this paper.

On the other hand, in the 1960s, Auslander \cite{Auslander} introduced a homological invariant for modules which is called Gorenstein dimension, or G-dimension for short.
After that, he further developed the theory of G-dimension with Bridger \cite{AB}.
Modules of G-dimension zero are called totally reflexive.
Over a Gorenstein local ring, totally reflexive modules are the same as maximal Cohen-Macaulay modules.

It is known that under a few assumptions, Gorenstein local rings of finite Cohen-Macaulay type are hypersurfaces \cite{Herzog}, and such rings and all nonisomorphic indecomposable maximal Cohen-Macaulay modules over them have been classified completely \cite{BGS}, \cite{GK}, \cite{Knorrer}.
(See also \cite{Yoshino}.)
Hence it is natural to ask whether there exists a non-Gorenstein local ring having only a finite number of isomorphism classes of indecomposable totally reflexive modules.
If such a ring exists, we want to determine all isomorphism classes of indecomposable totally reflexive modules.
However, the author conjectures that such a ring cannot essentially exist:

\begin{conj}\label{}
Let $R$ be a non-Gorenstein local ring.
Suppose that there is a nonfree totally reflexive $R$-module.
Then there are infinitely many isomorphism classes of indecomposable totally reflexive $R$-modules.
\end{conj}

Here, the assumption of the existence of a nonfree totally reflexive module is essential because, for example, all totally reflexive modules are free over a Cohen-Macaulay non-Gorenstein local ring with minimal multiplicity \cite{Yoshino2}.
The author proved that the above conjecture is true for any henselian local ring of depth at most two \cite{Takahashi0}, \cite{Takahashi1}, \cite{Takahashi2}.
He also proved recently that the conjecture is true for any henselian local ring having a nonfree cyclic totally reflexive module \cite{Takahashi3}.

In this paper, investigating the relationship between the number of totally reflexive modules and the dimensions of prime ideals, we consider the number of isomorphism classes of indecomposable totally reflexive modules.
The following result is (part of) the main theorem of this paper:

\begin{thm*}\label{}
Let $R$ be a local ring which is either complete or has uncountable residue field.
Suppose that there is a prime ideal $\p$ of $R$ with $\grade\,\p >0$ and $\dim R/\p >1$, and a totally reflexive $R$-module $M$ such that $M_\p$ is not $R_\p$-free.
Then there are uncountably many nonisomorphic indecomposable totally reflexive $R$-modules.
\end{thm*}

In the next section, we will show Theorem \ref{hl} without the assumption that the base ring $R$ is excellent.
In the last section, we shall prove our main theorem, and show that it implies Theorem \ref{hl} for a Cohen-Macaulay local ring with a canonical module.
We will also produce examples of rings over which there are an uncountably infinite number of nonisomorphic totally reflexive modules.

\section{On the theorem of Huneke and Leuschke}

Throughout this section, let $R$ be a commutative and noetherian ring, and all $R$-modules are considered to be finitely generated.

In this section, we consider Theorem \ref{hl}, which is due to Huneke and Leuschke.
Refining their proof, we check that the theorem holds without assuming that the base ring is excellent.
As Huneke and Leuschke do, we also need the following lemma, which is what is called ``countable prime avoidance''.
For the details, we refer to \cite[Corollaries (2.2),(2.6)]{SV}.
See also \cite[Lemma 3]{Burch}.

\begin{lem}[Countable Prime Avoidance]\label{avoid}
Let $(R, \m , k)$ be a local ring, and assume that $R$ is complete or that $k$ is an uncountable set.
Let $I$ be an ideal of $R$ and $\{\p_i\} _{i\in\N}$ a family of prime ideals of $R$.
If $I \nsubseteq\p _i$ for any $i\in\N$, then $I\nsubseteq\bigcup _{i\in\N}\p _i$.
\end{lem}

Recall that a subset $T$ of $\Spec\,R$ is said to be {\it stable under specialization} provided that if $\p\in T$ and $q\in\Spec\,R$ with $\p\subseteq\q$ then $\q\in T$.
For a ring satisfying countable prime avoidance, any prime ideal in a countable set which is stable under specialization has dimension at most one:

\begin{lem}\label{stsp}
Let $(R,\m,k)$ be a local ring, and assume either that $R$ is complete or that $k$ is uncountable.
Let $T$ be a countable subset of $\Spec\,R$ which is stable under specialization.
Then $\dim R/\p\leq 1$ for any $\p\in T$.
\end{lem}

\begin{pf}
Since the maximal ideal $\m$ is not contained in every prime ideal $\p\in T-\{\m\}$ and $T-\{\m\}$ is a countable set, Lemma \ref{avoid} implies that $\m$ is not contained in $\bigcup _{\p\in T-\{\m\}}\p$.
Hence we can choose an element $f\in\m$ which is not contained in any $\p\in T-\{\m\}$.

We claim that the ideal $\p +(f)$ is $\m$-primary for any $\p\in T$.
Indeed, take a prime ideal $\q$ containing $\p +(f)$.
As $\p$ belongs to $T$ and is contained in $\q$, the prime ideal $\q$ also belongs to $T$ by the assumption that $T$ is stable under specialization.
Since $f\in\q$, the choice of $f$ shows that $\q$ must coincide with $\m$.
Thus $\p +(f)$ is an $\m$-primary ideal.
We have $0=\dim (R/\p +(f))\geq\dim R/\p -1$, hence $\dim R/\p\leq 1$, which completes the proof of the lemma.
\qed
\end{pf}

We denote by $\mod\,R$ the category of finitely generated $R$-modules.
For a full subcategory $\X$ of $\mod\,R$, we denote by $\ind\,\X$ the set of isomorphism classes of $R$-modules in $\X$.
For a Cohen-Macaulay local ring $R$, set $\CM (R)$ to be the full subcategory of $\mod\,R$ consisting of all maximal Cohen-Macaulay $R$-modules.

Let $\Omega _R ^n$ be the $n$th syzygy functor over $R$.
(It is only well-defined up to free summands, but this is no restriction in the rest of the paper.)
Put $\Omega _R = \Omega _R^1$.
Every prime ideal in the singular locus of a Cohen-Macaulay local ring is determined by two indecomposable maximal Cohen-Macaulay modules.

\begin{prop}\label{sing}
Let $(R,\m,k)$ be a Cohen-Macaulay local ring.
Then
$$
\Sing\,R\subseteq\{\,\Ann\,\Tor _1 (X,Y)\,|\,X,Y\in\ind\,\CM(R)\,\}.
$$
\end{prop}

\begin{pf}
Fix a prime ideal $\p\in\Sing\,R$.
Set $M=\Omega _R ^d (R/\p )$ where $d=\dim R$.
Note that $M$ is a maximal Cohen-Macaulay $R$-module.

We claim that $\p = \Ann\,\Tor _1(M,M)$.
In fact, the module $\Tor _1(M,M)$ is isomorphic to $\Tor _{1+2d}(R/\p , R/\p )$, which is annihilated by $\p$.
Hence $\p\subseteq\Ann\,\Tor _1(M,M)$.
On the other hand, there are isomorphisms $\Tor _1(M,M)_\p\cong\Tor _{1+2d}(R/\p , R/\p )_\p\cong\Tor _{1+2d}^{R_\p}(\kappa (\p), \kappa (\p))$, where $\kappa (\p)$ denotes the residue field of $R_\p$.
Since the local ring $R_\p$ is not regular, the $R_\p$-module $\kappa (\p)$ has infinite projective dimension, and hence $\Tor _{1+2d}^{R_\p}(\kappa (\p), \kappa (\p))\neq 0$.
It follows from this that $\p$ belongs to $\Supp\,\Tor _1(M,M)$, therefore $\p\supseteq\Ann\,\Tor _1(M,M)$.
Thus we obtain $\p = \Ann\,\Tor _1(M,M)$.

Let $M=M_1\oplus M_2\oplus\cdots\oplus M_n$ be an indecomposable decomposition.
Each $M_i$ is an indecomposable maximal Cohen-Macaulay $R$-module.
Then $\p =\Ann\,\bigoplus _{1\leq i, j\leq n}\Tor _1(M_i, M_j)=\bigcap _{1\leq i, j\leq n}\Ann\,\Tor _1(M_i,M_j)$.
Noting that $\p$ is a prime ideal, we easily see that $\p =\Ann\,\Tor _1(M_a,M_b)$ for some integers $a, b$.
\qed
\end{pf}

Now, let us observe Theorem \ref{hl}.
Let $R$ be a Cohen-Macaulay local ring which is either complete or has uncountable residue field, and suppose that $R$ is of countable Cohen-Macaulay type.
Then $\ind\,\CM (R)$ is a countable set, and Proposition \ref{sing} implies that so is the singular locus $\Sing\,R$.
Note that $\Sing\,R$ is stable under specialization.
The theorem of Huneke and Leuschke is obtained from Lemma \ref{stsp} without the assumption that $R$ is excellent:

\begin{thm}\label{sch}
Let $(R, \m , k)$ be a Cohen-Macaulay local ring of countable Cohen-Macaulay type.
Suppose either that $R$ is complete or that $k$ is uncountable.
Then $\dim R/\p\leq 1$ for any $\p\in\Sing\,R$.
\end{thm}

\begin{rem}\label{partition}
Let $R$ be an arbitrary ring.
Let $\sf A$ be a set of nonisomorphic $R$-modules, and $\sf B$ the set of finite direct sums of $R$-modules in $\sf A$.
Then it is elementarily seen that the set $\sf A$ is countable if and only if so is $\sf B$.
Thus, in particular, the set of isomorphism classes of objects of a full subcategory $\X$ of $\mod\,R$ is countable if and only if so is the set of isomorphism classes of {\em indecomposable} objects of $\X$.
Hence, for the purpose of proving Theorem \ref{sch}, the conclusion of Proposition \ref{sing} is itself not necessary; Theorem \ref{sch} actually follows from the statement
$$
\Sing\,R\subseteq\{\,\Ann\,\Tor _1 (M,M)\,|\,M\in\CM(R)\,\},
$$
which we obtained in the middle of the proof of Proposition \ref{sing}.
\end{rem}

\section{Main theorem}

Throughout this section, $R$ is always assumed to be a commutative noetherian ring.
We also assume that all $R$-modules are finitely generated.

In this section we observe Theorem \ref{sch} from a more general viewpoint.
To be precise, we shall consider countablity of the set of isomorphism classes of indecomposable totally reflexive modules over an arbitrary local ring, and prove an analogue of Theorem \ref{sch}.
We begin with recalling the definitions of a semidualizing module and a totally reflexive module.

\begin{defn}
(1) An $R$-module $C$ is called {\it semidualizing} if the natural homomorphism $R\to\Hom _R(C,C)$ is an isomorphism and $\Ext _R^i(C,C)=0$ for any $i>0$.\\
(2) Let $C$ be a semidualizing $R$-module.
We say that an $R$-module $M$ is {\it totally $C$-reflexive} (or $M$ has {\it G$_C$-dimension zero}) if the natural homomorphism $M \to \Hom _R(\Hom _R(M,C),C)$ is an isomorphism and $\Ext _R ^i(M,C)=\Ext _R ^i(\Hom _R (M,C) ,C)=0$ for any $i>0$.
\end{defn}

For a semidualizing $R$-module $C$, we denote by $\G _C(R)$ the full subcategory of $\mod\,R$ consisting of all totally $C$-reflexive $R$-modules, and set
$$
\begin{cases}
W_C(R)=\{\,\p\in\Spec\,R\,|\,M_\p\text{ is not }R_\p\text{-free for some }M\in\G _C(R)\,\},\\
W_C^0(R)=\{\,\p\in W_C(R)\,|\,\grade\,\p >0\,\}.
\end{cases}
$$
Here we check several basic properties.

\begin{prop}\label{sd}
\begin{enumerate}
\item[{\rm (1)}]
A free $R$-module of rank one is semidualizing.
\item[{\rm (2)}]
Any free $R$-module is totally $C$-reflexive for every semidualizing $R$-module $C$.
\item[{\rm (3)}]
Let $R$ be a Cohen-Macaulay local ring with a canonical module $K$.
Then
\begin{enumerate}
\item[{\rm (i)}]
$K$ is a semidualizing $R$-module.
\item[{\rm (ii)}]
An $R$-module is totally $K$-reflexive if and only if it is maximal Cohen-Macaulay.
In other words, one has $\G_K(R)=\CM(R)$.
\item[{\rm (iii)}]
One has $W_K(R)=\Sing\,R$.
\end{enumerate}
\end{enumerate}
\end{prop}

\begin{pf}
The first and second statements immediately follow from definition.
As for the third statement, the assertions (i) and (ii) are basic properties of a canonical module; see \cite[Theorem 3.3.10]{BH} for example.

Let us prove the assertion (iii).
Fix a prime ideal $\p$ of $R$.
If $\p$ is in $W_K(R)$, then $M_\p$ is not $R_\p$-free for some maximal Cohen-Macaulay $R$-module $M$ by (ii).
Note that $M_\p$ is maximal Cohen-Macaulay over $R_\p$, and that any maximal Cohen-Macaulay module over a regular local ring is free.
Hence $\p$ belongs to $\Sing\,R$.
Conversely, suppose that this is the case.
Putting $M=\Omega _R ^d (R/\p)$ where $d=\dim R$, we see that $M$ is a maximal Cohen-Macaulay $R$-module, and that $M_\p$ is isomorphic to the $d$th syzygy of the resudue field $\kappa (\p)$ of $R_\p$ up to free summand.
Since the local ring $R_\p$ is not regular, no syzygy of $\kappa (\p)$ is free.
Therefore $M_\p$ is not a free $R_\p$-module, that is, the prime ideal $\p$ belongs to $W _K (R)$.
\qed
\end{pf}

A totally $R$-reflexive $R$-module is simply called {\it totally reflexive} (or {\it G-dimension zero}).
We put $\G (R)=\G _R(R)$, $W(R)=W_R(R)$ and $W^0(R)=W_R^0(R)$.

For an $R$-module $M$, let $\NF (M)$ denote the {\it nonfree locus} of $M$:
$$
\NF (M) = \{\,\p\in\Spec\,R\,|\,M_\p\text{ is not }R_\p\text{-free}\,\}.
$$
This set is closed in $\Spec\,R$, which is well-known to experts.

\begin{lem}\label{closed}
For any $R$-module $M$, one has
$$
\NF (M) = \Supp\,\Ext ^1(M,\Omega M) = V(\Ann\,\Ext ^1(M,\Omega M)).
$$
In particular, $\NF (M)$ is a closed subset of $\Spec\,R$.
\end{lem}

\begin{pf}
Let $\p\in\Supp\,\Ext ^1(M,\Omega M)$.
Then $\Ext _{R_\p}^1(M_\p,(\Omega M)_\p )\cong (\Ext _R^1(M, \Omega M))_\p\neq 0$, and $M_\p$ is not $R_\p$-free.
Conversely, let $\p\in\NF (M)$.
There is an exact sequence $\sigma : 0 \to \Omega M \to F \to M \to 0$ of $R$-modules, where $F$ is free.
Since $M_\p$ is not $R_\p$-free, the localization $\sigma _\p : 0 \to (\Omega M)_\p \to F_\p \to M_\p \to 0$ does not split.
Hence $\sigma _\p$ corresponds to a nonzero element of $\Ext _{R_\p}^1(M_\p,(\Omega M)_\p )$, and we have $\Ext _{R_\p}^1(M_\p,(\Omega M)_\p )\neq 0$.
Therefore $\p$ is in $\Supp\,\Ext ^1(M,\Omega M)$.
\qed
\end{pf}

The following result will play a key role in the rest of this paper.

\begin{lem}\label{key}
For any $\p\in W_C^0(R)$ there exists $M\in\G _C(R)$ such that
$$
\NF (M) = V(\p ).
$$
\end{lem}

\begin{pf}
Fix a prime ideal $\p$ in $W_C^0(R)$.
By definition there exists $M\in\G _C(R)$ such that $M_\p$ is not $R_\p$-free, i.e., $\p\in\NF (M)$.
Noting that $\NF (M)$ is stable under specialization, we see that $\NF (M)$ contains $V(\p )$.
If $\NF (M)$ coincides with $V (\p )$, then there is nothing to prove.
Hence let us assume that $\NF(M)$ strictly contains $V(\p)$, and take a prime ideal $\q\in\NF(M)-V(\p)$.
The prime ideal $\p$ is contained neither in $\q$ nor in any $\P\in\Ass\,R$ since $\p$ has positive grade.
Hence we can choose an $R$-regular element $x$ in $\p -\q$.
The isomorphism $R\cong\Hom (C,C)$ shows that $\Ass\,R=\Ass\,C$, which implies that the element $x$ is $C$-regular.
Hence there is an exact sequence $0 \to C \overset{x}{\to} C$, and applying the functor $\Hom (\Hom (M,C), -)$ to this, we see that $x$ is also $M$-regular.
Put $N=\Omega _R(M/xM)$.

\setcounter{claim}{0}
\begin{claim}\label{ex}
The $R$-module $N$ belongs to $\G _C(R)$.
\end{claim}

\begin{cpf}
We have two exact sequences $0 \to M \overset{x}{\to} M \to M/xM \to 0$ and $0 \to N \to F \to M/xM \to 0$, where $F$ is a free $R$-module.
From these sequences we make the following pullback diagram:
$$
\begin{CD}
@. @. 0 @. 0 \\
@. @. @VVV @VVV \\
@. @. N @= N \\
@. @. @VVV @VVV \\
0 @>>> M @>>> P @>>> F @>>> 0 \\
@. @| @VVV @VVV \\
0 @>>> M @>{x}>> M @>>> M/xM @>>> 0 \\
@. @. @VVV @VVV \\
@. @. 0 @. 0
\end{CD}
$$
Since the middle row splits, $P$ is isomorphic to $M\oplus F$, and we get an exact sequence
$$
0 \to N \to M\oplus F \to M \to 0.
$$
As both $M$ and $F$ belong to $\G _C(R)$, so does $N$ by \cite[Theorem 2.1(1)]{ATY}.
\qed
\end{cpf}

\begin{claim}\label{ex2}
One has $\NF (M)\supsetneq\NF (N)\supseteq V(\p)$.
\end{claim}

\begin{cpf}
First of all, we notice that $N_\p$ (resp. $N_\q$) is isomorphic to $\Omega _{R_\p}(M_\p /xM_\p )$ (resp. $\Omega _{R_\q}(M_\q /xM_\q )$) up to $R_\p$-free (resp. $R_\q$-free) summand.
Since $x$ is not in $\q$, we have $xM_\q = M_\q$ and $N_\q$ is $R_\q$-free, i.e., $\q\notin\NF (N)$.
Recall that $\q\in\NF (M)$.
Hence, in particular, the set $\NF (M)$ does not coincide with $\NF (N)$.

Suppose that $N_\p$ is $R_\p$-free.
Then the $R_\p$-module $M_\p /xM_\p$ has projective dimension at most one.
Since $x$ is in $\p$ and is $M$-regular, $x$ is $M_\p$-regular as an element of $R_\p$.
Hence we have $\pd _{R_\p}\,M_\p /xM_\p = \pd _{R_\p}\,M_\p +1$, and see that $M_\p$ is $R_\p$-free.
This is a contradiction, which shows that $N_\p$ is not $R_\p$-free, namely, $\p\in\NF(N)$.
Thus we obtain $V(\p)\subseteq\NF(N)$.

It remains to prove that $\NF (M)$ contains $\NF (N)$.
Let $\P\in\NF(N)$.
Then $\Omega _{R_\P}(M_\P /xM_\P )$ is a nonfree $R_\P$-module.
Hence $x$ must belong to $\P$.
Assume that $\P$ does not belong to $\NF (M)$.
Then $M_\P$ is a free $R_\P$-module.
Noting that $x$ is an $R_\P$-regular element, we easily see that $\Omega _{R_\P}(M_\P /xM_\P )$ is $R_\P$-free, and get a contradiction.
It follows that $\P$ belongs to $\NF (M)$, as desired.
\qed
\end{cpf}

If $\NF(N)$ strictly contains $V(\p)$, then in a similar way as above we can obtain $N'\in\G _C(R)$ such that $\NF(N)\supsetneq\NF(N')\supseteq V(\p)$.
Iterating this procedure yields a strict descending chain
$$
\NF (M)\supsetneq\NF (N)\supsetneq\NF (N')\supsetneq\cdots
$$
of closed subsets of $\Spec\,R$ (see Lemma \ref{closed}).
Since $\Spec\,R$ is a noetherian space, we cannot iterate the procedure infinitely many times.
Therefore there exists $L\in\G _C(R)$ such that $\NF (L)$ coincides with $V (\p )$.
\qed
\end{pf}

Making use of Lemma \ref{key}, we can represent each prime ideal in $W_C^0(R)$ by two indecomposable totally $C$-reflexive $R$-modules.

\begin{prop}\label{w}
$W_C^0(R)\subseteq\left\{\,\sqrt{\Ann\,\Ext ^1(X,Y)}\,\Biggm\vert\,X,Y\in\ind\,\G _C(R)\,\right\}$.
\end{prop}

\begin{pf}
Let $\p\in W_C^0(R)$.
According to Lemma \ref{key}, there is an $R$-module $M$ in $\G_C(R)$ with $\NF(M)=V(\p)$.
By Lemma \ref{closed} we have $\p = \sqrt{\Ann\,\Ext ^1(M,\Omega M)}$.
Decompose $M$ and $\Omega M$ into indecomposable $R$-modules:
$$
\begin{cases}
M = X_1\oplus\cdots\oplus X_n, \\
\Omega M = Y_1\oplus\cdots\oplus Y_m.
\end{cases}
$$
Here each $X_i$ and each $Y_j$ belong to $\ind\,\G_C(R)$.
Hence
$$
\p = \sqrt{\Ann\,\bigoplus _{1\leq i\leq n,1\leq j\leq m}\Ext ^1(X_i,Y_j)}=\bigcap _{1\leq i\leq n,1\leq j\leq m}\sqrt{\Ann\,\Ext ^1(X_i,Y_j)}.
$$
Since $\p$ is a prime ideal, we conclude that $\p = \sqrt{\Ann\,\Ext ^1(X_a,Y_b)}$ for some integers $a,b$.
This proves the proposition.
\qed
\end{pf}

Now, we have reached the stage to accomplish our aim in this section; the following theorem is our main result.

\begin{thm}\label{main}
Let $(R, \m , k)$ be a local ring, and assume either that $R$ is complete or that $k$ is uncountable.
Let $C$ be a semidualizing $R$-module such that $\ind\,\G _C (R)$ is a countable set.
Then the following hold.
\begin{enumerate}
\item[{\rm (1)}]
One has $\dim R/\p\leq 1$ for any $\p\in W_C^0(R)$.
\item[{\rm (2)}]
If $R$ satisfies Serre's $(S_1)$-condition, one has $\dim R/\p\leq 1$ for any $\p\in W_C (R)$.
\end{enumerate}
\end{thm}

\begin{pf}
First of all, note that both $W_C(R)$ and $W_C^0(R)$ are stable under specialization.
Since $\ind\,\G _C (R)$ is countable, Proposition \ref{w} shows that so is the set $W_C^0(R)$.
It follows from Lemma \ref{stsp} that $\dim R/\p\leq 1$ for any $\p\in W_C^0(R)$.
This proves the first assertion of the theorem.
As to the second assertion, notice that $R$ satisfies $(S_1)$ if and only if all the prime ideals of $R$ of grade zero are minimal primes.
Hence, if this is the case, then $W_C(R)$ is contained in the union set $W_C^0(R)\cup\Min\,R$, which is countable because the set $\Min\,R$ is finite.
Thus $W_C(R)$ is also countable, and Lemma \ref{stsp} implies that $\dim R/\p\leq 1$ for any $\p\in W_C(R)$.
\qed
\end{pf}

The above theorem shows Theorem \ref{sch} in the case where the base ring admits a canonical module:

\begin{cor}\label{}
Let $(R, \m , k)$ be a Cohen-Macaulay local ring of countable Cohen-Macaulay type with a canonical module $K$.
Suppose either that $R$ is complete or that $k$ is uncountable.
Then $\dim R/\p\leq 1$ for any $\p\in\Sing\,R$.
\end{cor}

\begin{pf}
Since $R$ is Cohen-Macaulay, $R$ satisfies $(S_1)$.
It is seen by the assumption and Proposition \ref{sd}(3)(ii) that $\ind\,\G_K(R)$ is a countable set.
By virtue of Theorem \ref{main}(2) and Proposition \ref{sd}(3)(iii), we get the conclusion of the corollary.
\qed
\end{pf}

\section{Examples}

We end this paper by giving several examples of rings having an uncountably infinite number of nonisomorphic indecomposable totally reflexive modules.

\begin{ex}\label{ex1}
Let
$$
R=k[[x,y,z,w]]/(x^2,yz,yw)
$$
where $k$ is an arbitrary field.
Then $R$ is a complete local ring of dimension two and depth one.
In particular, $R$ is not Cohen-Macaulay.
Since $(x^2,yz,yw)=(x^2,y)\cap (x^2,z,w)$ in $k[[x,y,z,w]]$, the associated primes of $R$ are $\p=(x,y)R$ and $\q=(x,z,w)R$, both of which are minimal.
Hence $\Ass\,R=\Min\,R$, equivalently, $R$ satisfies $(S_1)$.
There is an exact sequence
$$
\cdots \overset{x}{\to} R \overset{x}{\to} R \overset{x}{\to} \cdots
$$
whose $R$-dual sequence is exact.
Therefore $M:=R/xR$ is a totally reflexive $R$-module.
It is easy to see that $M_\p$ is not $R_\p$-free, hence $\p\in W(R)$.
Since $\dim R/\p =2$, we see by the second assertion of Theorem \ref{main} that $\ind\,\G(R)$ is an uncountable set.
\end{ex}

\begin{ex}\label{ex2}
Let $(S,\n)$ be a complete local ring of positive depth, and let
$$
R=S[[x,y,z]]/(x^2).
$$
Then the Gorenstein property or the Cohen-Macaulay property of $R$ is equivalent to that of $S$.
Set $\p = \n R+xR$.
This is a prime ideal of $R$ and $\dim R/\p =2$.
Since $S$ has positive depth, there is an $S$-regular element $f\in\n$.
Noting that $R$ is faithfully flat over $S$, we see that $f$ is also $R$-regular and $\p$ has positive grade because it contains $f$.
Put $M=R/xR$.
We see that $M$ is a totally reflexive $R$-module and that $M_\p$ is a nonfree $R_\p$-module.
Hence $\p\in W^0(R)$, and therefore it follows from the first assertion of Theorem \ref{main} that $\ind\,\G(R)$ is an uncountable set.
\end{ex}

\begin{ex}\label{ex3}
Let $S$ be a complete local domain which is not a field, and let
$$
R=S[[x,y,z]]/(x^2).
$$
Then Example \ref{ex2} guarantees that there exist an uncountably infinite number of isomorphism classes of indecomposable totally reflexive $R$-modules.
In the following, let us actually construct such modules.

Set $M=R/xR$ and $\p = (x,y)R$.
The $R$-module $M$ is totally reflexive and $\p$ is a prime ideal of $R$ such that $\dim R/\p =\dim S+1\geq 2$.
Since $y$ is an $R$-regular element, $\p$ has positive grade.
The $R_\p$-module $M_\p$ is nonfree.
For an element $f\in S[[z]]\subseteq R$, put $\p ^f = (x, y-zf)R$.
Then $\p ^f$ is a prime ideal of $R$ and $\p ^f = \Omega _R(M/(y-zf)M)$.
Since $y-zf$ is an $R$-regular element, the proof of Claim \ref{ex} in the proof of Lemma \ref{key} shows that $\p ^f$ is a totally reflexive $R$-module.

\setcounter{claim}{0}
\begin{claim}\label{623}
$\p ^f$ is indecomposable as an $R$-module.
\end{claim}

\begin{cpf}
Set $T=S[[x,y,z]]$ and denote by $\M$ the maximal ideal of $T$.
Assume that $\p ^f$ is decomposable, and write $\p ^f=\alpha R\oplus\beta R$ for some $\alpha , \beta \in T$.
Then we have $(x, y-zf)R=(\alpha, \beta)R$.
Denote by $\overline{t}$ the residue class of $t\in T$ in $R=T/x^2T$.
The elements $\overline{x}, \overline{y-zf}$ are part of a minimal system of generators of the maximal ideal $\m$ of $R$, hence so are $\overline{\alpha},\overline{\beta}$.
Therefore $\overline{\alpha}, \overline{\beta}\in\m -\m ^2$, which implies that $\alpha,\beta\in\M -\M ^2$.
We have $\overline{\alpha}\overline{\beta}\in\alpha R\cap\beta R=0$, namely $\alpha\beta\in x^2T$.
Noting that $xT$ is a prime ideal of $T$, we get either $\alpha\in xT$ or $\beta\in xT$.
We may assume $\alpha\in xT$, and write $\alpha = x\gamma$ for some $\gamma\in T$.
Since $x\in\M$ and $\alpha\notin\M ^2$, the element $\gamma$ is a unit of $T$, and we obtain $\beta \in xT$ since $\alpha\beta\in x^2T$.
However, we then have $y-zf\in\p ^f=(\alpha,\beta)R\subseteq xR$, which is a contradiction.
\qed
\end{cpf}

\begin{claim}\label{642}
$\NF (\p ^f) = V(\p ^f)$.
\end{claim}

\begin{cpf}
For any prime ideal $\q$ not containing $\p ^f$, one has $(\p ^f)_\q = R_\q$.
This shows that $V(\p ^f)$ contains $\NF (\p ^f)$.
On the other hand, we easily see that $R_{\p ^f}$ is not a regular local ring.
Hence $(\p ^f)_{\p ^f}=\p ^fR_{\p ^f}$ is not $R_{\p ^f}$-free, which implies that $\p ^f$ belongs to $\NF(\p ^f)$.
Since $\NF (\p ^f)$ is stable under specialization, the set $\NF(\p ^f)$ contains $V(\p ^f)$.
\qed
\end{cpf}

\begin{claim}\label{655}
The condition $\p ^f=\p ^g$ implies $f=g$.
\end{claim}

\begin{cpf}
For an element $h\in S[[z]]$, we define a homomorphism
$$
\phi ^h: R/\p ^h \to S[[z]]
$$
by $\phi ^h(\overline y)=zh$, $\phi ^h(\overline z)=z$ and $\phi ^h(\overline s)=s$ for any $s\in S$.
Here, for $r \in R$, $\overline r$ denotes the residue class of $r$ in $R/\p ^f$.
Then it is easily seen that $\phi ^h$ is an isomorphism.

Now, suppose that $\p ^f =\p ^g$.
Let $\Phi : S[[z]] \to S[[z]]$ be the composite map $\phi ^g\cdot (\phi ^f)^{-1}$.
Then we notice that this map $\Phi$ is an identity map since it sends $z$ and each $s\in S$ to themselves.
Noting that $\Phi (zf)=zg$, we see that $f=g$.
\qed
\end{cpf}

It follows from Claims \ref{642} and \ref{655} that for $f,g\in S$ one has $f=g$ if and only if $\p ^f=\p ^g$, if and only if $V(\p ^f)=V(\p ^g)$, if and only if $\NF (\p ^f)=\NF (\p ^g)$.
Hence we see by Claim \ref{623} that for any two distinct elements $f,g\in S$ the prime ideals $\p ^f$ and $\p ^g$ of $R$ are indecomposable totally reflexive $R$-modules which are not isomorphic to each other.
In other words, the map from $S[[z]]$ to the set of nonisomorphic indecomposable totally reflexive $R$-modules which sends $f$ to $\p ^f$ is injective.
Since $S[[z]]$ is an uncountable set, so is the set $\{ \p ^f \}_{f\in S[[z]]}$.
Thus we obtain uncountably many nonisomorphic totally reflexive $R$-modules.
\end{ex}


{\sc Acknowledgments.}
The author expresses his gratitude to the participants in the Seminar on Commutative Ring Theory at Meiji University, which is organized by Shiro Goto.
Thanks to their comments and suggestions, this paper is in the present simple form.


\end{document}